\documentclass[12pt,a4paper]{article}
\usepackage{amsmath,amsthm}
\usepackage{amssymb}

\usepackage{enumitem}

\usepackage{graphicx}

\usepackage[T1]{fontenc}
\newtheorem{theorem}{Theorem}[section]
\newtheorem{corollary}[theorem]{Corollary}
\newtheorem{lemma}[theorem]{Lemma}
\newtheorem{proposition}[theorem]{Proposition}

\theoremstyle{definition}
\newtheorem{definition}[theorem]{Definition}
\newtheorem{remark}[theorem]{Remark}
\newtheorem{example}[theorem]{Example}

\numberwithin{equation}{section}

\DeclareMathOperator{\Ped}{Ped}
\DeclareMathOperator{\her}{her}
\DeclareMathOperator{\id}{id}

\begin{document}

%%%%% To ease editing, for IMPAN journals add:

\baselineskip=17pt

%%%%%%%%%%%%%%%%

\title{Remarks on some simple $C^*$-algebras admitting a unique lower semicontinuous $2$-quasitrace}

\author{Jacopo Bassi\\
Department of Mathematics\\ 
University of Tor Vergata\\
Via della Ricerca Scientifica 1\\
00133 Roma, Italy\\
E-mail: bassi@mat.uniroma2.it}
%\and 
%Jan Krzysztof Nowak\\
%Institute of Mathematics\\ 
%{\L}ojasiewicza 6\\
%Jagiellonian University\\
%30-348 Krak\'ow, Poland\\
%E-mail: jk.nowak@im.uj.edu.pl}

\date{}

\maketitle

%% Classification and key words; note that the 2010 classification is used:

\renewcommand{\thefootnote}{}



\renewcommand{\thefootnote}{\arabic{footnote}}
\setcounter{footnote}{0}

%%%%%%%%

\begin{abstract}
Using different descriptions of the Cuntz semigroup and of the Pedersen ideal, it is shown that $\sigma$-unital simple $C^*$-algebras with almost unperforated Cuntz semigroup, a unique lower semicontinuous $2$-quasitrace and whose stabilization has almost stable rank $1$ are either stable or algebraically simple.
\end{abstract}

\section{Introduction}  \label{sect1}
The Cuntz semigroup was introduced by Cuntz in \cite{cuntzdim} in order to study the existence of quasitraces on certain $C^*$-algebras. Since then, it has received attention in view of the couterexamples to the classification program of $C^*$-algebras by means of their Elliott invariant constructed by Toms in \cite{toms}, where it is proved the existence of simple, unital, nuclear, non-elementary $C^*$-algebras sharing the same Elliott invariant, but whose Cuntz semigroups differ. Despite the existence of many non-isomorphic nuclear, simple $C^*$-algebras with the same Cuntz semigroup, in \cite{nccw} it was proved that a modified version of the Cuntz semigroup classifies certain inductive limits of $1$-dimensional $NCCW$ (Non-Commutative CW) complexes.\\

In \cite{cuntz_hm} it is shown that the Cuntz semigroup of a $C^*$-algebra can be equivalently described in terms of certain classes of Hilbert modules, and in \cite{open_proj} another characterization in terms of classes of open projections in the double dual of the $C^*$-algebra is given. In both these papers it is observed that, under the assumption of almost stable rank $1$, Cuntz comparison can be reduced to a simpler relation. Also, as shown in \cite{lowersemi}, \cite{cuntz_nonunital} and \cite{thiel_rank}, in some cases, the Cuntz semigroup of a $C^*$-algebra can be reconstructed from the lower semicontinuous functions on the cone of lower semicontinuous ($2$-quasi)traces on the $C^*$-algebra and the Murray-von Neumann semigroup.\\

The Cuntz semigroup was recently used in \cite{r1} and \cite{thiel_rank} in order to give answers to important open problems in the field of $C^*$-algebras, such as a conjecture by Blackadar and Handelman on dimension functions, the Global Glimm Halving Problem, and questions on the equivalence of certain different comparability properties.\\

Specializing to $C^*$-algebras with suitable properties, in the present article it is given a description of the Pedersen ideal in terms of Cuntz comparison, and the above correspondence between different descriptions of the Cuntz semigroup is employed in order to establish a connection between algebraic simplicity/stability and related finiteness/infiniteness properties.

\subsection{Notation} If $A$ is a $C^*$-algebra, we denote by $\tilde{A}$ its minimal unitization; if $A$ is unital, we denote by $GL(A)$ the group of invertible elements in $A$. The $C^*$-algebra of compact operators on a separable infinite-dimensional Hilbert space is denoted by $\mathbb{K}$.\\
 If $A$ is a $C^*$-algebra and $x$ is a positive element in $A$, we denote by $E_x=\overline{xA}$ the corresponding right ideal, viewed as a Hilbert module. If $x$ is strictly positive for $A$, we identify $E_x$ with $A$. The hereditary $C^*$-algebra generated by an element $y \in A$ is denoted by $\her (y)$.

\section{The Cuntz semigroup}
This section contains definitions and preliminary results needed in the forthcoming part. Let $A$ be a $C^*$-algebra and $a$, $b$ two positive elements in $A$; then \textit{$a$ is Cuntz subequivalent to $b$} ($a \precsim b$) if there is a sequence $\{x_n\}$ of elements in $A$ such that $a=\lim x_n b x_n^*$; and \textit{$a$ is Cuntz equivalent to $b$} ($a \sim b$) if $a \precsim b$ and $b \precsim a$. If $A$ is a $C^*$-algebra, its \textit{Cuntz semigroup} $(Cu (A), \oplus, \leq) = Cu (A)$ is the partially ordered abelian semigroup whose elements are Cuntz classes of positive elements in $A\otimes \mathbb{K}$, where the semigroup structure is obtained by composing any $*$-isomorphism $M_2 (\mathbb{K})\otimes A \rightarrow \mathbb{K} \otimes A$ with the orthogonal sum $ \mathbb{K} \otimes A \times \mathbb{K} \otimes A \rightarrow M_2 (\mathbb{K})\otimes A$. The partial order is induced by Cunts subequivalence: if $[a], [b] \in Cu (A)$ are the Cuntz classes of elements $a$ and $b$ respectively, then $[a] \leq [b]$ if and only if $a \precsim b$. This is the same as $W(A\otimes \mathbb{K})$ in the notation of \cite{cuntz-semigroup}.\\ 

If $A$ is a $C^*$-algebra, every increasing sequence in $Cu(A)$ admits a supremum (\cite{cuntz-semigroup} Proposition 4.17). Let $[a]$ and $[b]$ be elements in $Cu (A)$, we say that \textit{$[a]$ is compactly contained in $[b]$} ($[a] \ll [b]$) if for every increasing sequence $\{ [b_n]\} \subset Cu (A)$ such that $[b] \leq \sup_n [b_n]$, there is a $k \in \mathbb{N}$ such that $[a] \leq [b_k]$ (this is the way below relation in \cite{cuntz-semigroup} Definition 4.1).

\begin{definition}[\cite{cuntz_t} Section 2.1 and Definition 5.3.1] Let $A$ be a $C^*$-algebra and $[a] \in Cu(A)$. $[a]$ is \textit{compact} if $[a] \ll [a]$. $[a]$ is \textit{soft} if for every $[b] \in Cu(A)$ satisfying $[b] \ll [a]$ there is $n \in \mathbb{N}$ such that $(n+1) [b] \leq n [a]$.
\end{definition}

If $A$ is a $C^*$-algebra, the purely non-compact elements in $Cu(A)$, defined in \cite{cuntz_hm} pag. 27, are soft by \cite{cuntz_t} Proposition 5.3.5.

\begin{definition}[\cite{lowersemi}, Introduction] A \textit{$2$-quasitrace} on a $C^*$-algebra $A$ is a map $\tau: (A\otimes \mathbb{K})_+ \rightarrow [0,\infty]$ that sends $0$ to $0$, is linear on commuting elements and satisfies the trace identity: $\tau (x^* x) = \tau (xx^*)$ for every $x \in A \otimes \mathbb{K}$.
\end{definition}

Denote the set of lower semicontinuous $2$-quasitraces on $A$ by $QT_2 (A)$.

\begin{definition}[\cite{lowersemi} pag. 17]
A \textit{functional} on the Cuntz semigroup of a $C^*$-algebra $A$ is a map $Cu (A) \rightarrow [0,\infty]$ that preserves the $0$ element, sum, partial order, and suprema of increasing sequences.
\end{definition}
We denote the set of functionals on $Cu (A)$ by $F(Cu (A))$. As usual, if $\lambda$ is a functional on $Cu(A)$, we say that $\lambda$ is \textit{trivial} if it only takes values $0$ and $\infty$ and is \textit{nontrivial} otherwise. The same terminology is used for the lower semicontinuous $2$-quasitraces. There is a well known bijection $QT_2 (A) \simeq F(Cu(A))$ (\cite{lowersemi} Proposition 4.2).\\

If $A$ is a simple $C^*$-algebra admitting a non-trivial lower semicontinuous $2$-quasitrace, every non-zero element in $Cu(A)$ is either compact or soft by \cite{cuntz_t} Proposition 5.3.16.\\

Let $A$ be a $C^*$-algebra and consider the trivial lower semicontinuous $2$-quasitraces on $A$ defined by
\[
\tau_0 (a) = 0 \quad \mbox{ for every } a \in (A\otimes \mathbb{K})_+, \qquad
\tau_\infty (a) = \begin{cases}	0	&	\mbox{ for } a=0\\
						\infty	&	\mbox{ for } a \in (A\otimes \mathbb{K})_+ \backslash \{ 0 \}.
						\end{cases}.
\]
In the same way, we define trivial functionals on $Cu(A)$ by
\[
\lambda_0 ([a]) = 0 \quad \mbox{ for all } [a] \in Cu (A), \qquad
\lambda_\infty ([a]) = \begin{cases}	0	&	\mbox{ for } [a]=0 \\
							\infty	&	\mbox{ for } [a] \in Cu (A) \backslash \{ 0 \}
							\end{cases}.
							\]

\begin{lemma}
\label{nontrivial}
Let $A$ be a simple $C^*$-algebra. Then the only possible trivial lower semicontinuous $2$-quasitraces on $A$ are $\tau_0$ and $\tau_\infty$; the only possible trivial functionals on $Cu(A)$ are $\lambda_0$ and $\lambda_\infty$.\\
Furthermore, $A$ admits a unique nontrivial lower semicontinuous $2$-quasitrace (up to a scalar) if and only if $Cu (A)$ admits a unique nontrivial functional (up to a scalar).
\end{lemma}
\proof 
Let $\lambda$ be a trivial functional on $Cu(A)$. It follows from the proof of \cite{purinf} Lemma 3.12 that the set $N_\lambda := \{ a \in A\otimes \mathbb{K} \; : \; \lambda ([|a|]) =0 \}$ is an algebraic ideal. By simplicity it is either the set $\{ 0 \}$ or is dense in $A\otimes \mathbb{K}$.\\
If $N_\lambda = \{ 0\}$, then $\lambda = \lambda_\infty$; hence suppose $N_\lambda \neq0$. Let $[a]$ be an element of $Cu (A)$  and note that $\lambda ([a]) = \sup_{\epsilon >0} \lambda ([(a-\epsilon)_+])$. Since $N_\lambda$ is dense, $\Ped (A\otimes \mathbb{K}) \subset N_\lambda$ and so $\lambda([(a-\epsilon)_+])=0$ for every $\epsilon >0$. Thus $\lambda ([a])=0$ and $\lambda = \lambda_0$.\\
Using the bijective correspondence between functionals on the Cuntz semigroup and lower semicontinuous $2$-quasitraces on the corresponding $C^*$-algebra, the $2$-quasitrace associated to a functional $\lambda$ is $\tau_\lambda (a) = \int_0^\infty \lambda ([(a-t)_+])$ for $a \in A_+$. Hence the functional $\lambda_0$ corresponds to the lower semicontinuous $2$-quasitrace $\tau_0$ and $\lambda_\infty$ corresponds to $\tau_\infty$.\\
It follows that if $A$ admits a unique non-trivial lower semicontinuous $2$-quasitrace up to scalars, the multiples of it correspond to the same multiples of the associated non-trivial functional on $Cu(A)$. $\Box$\\

Another property of the Cuntz semigroup of a $C^*$-algebra which turns out to be important, in particular for its connection with the Toms-Winter conjecture, is the so called \textit{almost unperforation}:
\begin{definition}[\cite{lowersemi} pag. 26]
A partially ordered semigroup $(S,+, \leq )$ is said to be almost unperforated if, whenever there are elements $s$, $t$ in $S$ and $n \in \mathbb{N}$ such that $(n+1)s \leq nt$, then $s \leq t$.
\end{definition}
\begin{theorem}[\cite{cuntz_t} Theorem 5.3.12]
\label{comparison}
Let $A$ be a $C^*$-algebra such that $Cu(A)$ is almost unperforated. If $[a]$ and $[b]$ are elements in $Cu (A)$ such that $[a]$ is soft and $\lambda ([a]) \leq \lambda ([b])$ for every functional $\lambda \in F(Cu (A))$, then $[a] \leq [b]$.
\end{theorem}

As already mentioned in the Introduction, $C^*$-algebras with almost stable rank $1$ have received a particular attention in \cite{cuntz_hm} and \cite{open_proj}, since for these $C^*$-algebras Cuntz subequivalence can be described in terms of other relations. This fact will be an important ingredient in the following section. We recall this concept.
\begin{definition}[\cite{zstable_projless} Definition 3.1]
Let $A$ be a $C^*$-algebra. We say that $A$ has \textit{almost stable rank $1$} if for every hereditary $C^*$-subalgebra $B \subset A$, we have $B \subset \overline{GL(\tilde{B})}$.
\end{definition}

Note that almost stable rank $1$ is not a stable property (\cite{thiel_rank} Example 6.6). 

\section{Algebraic simplicity and stability}

If $A$ is a $C^*$-algebra, we denote by $\Ped (A)$ its Pedersen ideal, i.e. the dense ideal that is contained in every other dense ideal and refer to  \cite{pedersen} 5.6 for its construction and properties. We will give two descriptions of the Pedersen ideal for $C^*$-algebras belonging to a certain class and use these for the proof of Theorem \ref{dyc2}.\\
The first characterization of the Pedersen ideal is in terms of compact containment for a certain class of $C^*$-algebras. A more general result, that appears as a Remark on page 4 of \cite{cuntz_nonunital}, was suggested to hold by Elliott without the assumption that the stabilization of the $C^*$-algebra has almost stable rank 1.

\begin{proposition}
\label{ped}
Let $A$ be a $C^*$-algebra such that $A\otimes \mathbb{K}$ has almost stable rank $1$. Then
\[
\Ped (A) = \{ a \in A \; : \; \exists b \in M_\infty (A)_+ \; \mbox{ such that } \; [ | a | ] \ll [b] \}.
\]
\end{proposition}
\proof 
Let $a$ be an element of $A$ and write $a=u |a|$ for its polar decomposition in $A^{**}$. If $|a|$ belongs to $\Ped (A)$, then by \cite{pedersen} Proposition 5.6.2 $|a|^{\alpha}$ belongs to $\Ped (A)$ for every $\alpha >0$; note now that $u |a|^{1-\alpha}$ belongs to $A$ for every $0<\alpha <1$ and so $u|a| = u |a|^{1-\alpha} |a|^\alpha$ belongs to $\Ped (A)$. Thus we are left to prove that the set
\[
J_+ := \{ a \in A_+ \; : \; \exists b \in M_\infty (A)_+ \; \mbox{ such that } \; [  a  ] \ll [b] \\]
\]
is contained in $\Ped (A)$, that
\[
J := \{ a \in A \; : \; \exists b \in M_\infty (A)_+ \; \mbox{ such that } \; [ | a | ] \ll [b] \}
\]
is an ideal and that it is dense.\\
Let $a$ be a positive element in $A$ and suppose there is $n \in \mathbb{N}$ and a positive element $b$ in $M_n (A)$ such that $[a] \ll [b]$. Then there is $\epsilon >0$ such that $a \precsim (b-\epsilon)_+$ in $M_n (A)$. Since $A\otimes \mathbb{K}$ has almost stable rank $1$, by the proof of \cite{cuntz_hm} Theorem 3 and by Proposition 4.6 of \cite{open_proj}, there is an element $y$ in $M_n (A)$ such that $\her(a) = \her(yy^*)$, $y^* y \in \her ((b-\epsilon)_+)$; in particular by \cite{pedersen} Proposition 5.6.2 $y^*y$ belongs to $\Ped (M_n (A)) $ and the same is true for $yy^*$.\\
Now we want to show that for any $n \in \mathbb{N}$, $\Ped (M_n (A)) = \Ped (A) \otimes_{alg} M_n$. Notice that $\Ped (A) \otimes_{alg} M_n$ is a dense ideal in $M_n (A)$ and so $\Ped (M_n (A)) \subset \Ped (A) \otimes_{alg} M_n$. On the other hand, if $x = (x_{i,j})$ is an element in $M_n (A)$ such that $x_{i,j}$ belongs to $\Ped (A)$ for every $1 \leq i,j \leq n$, we can write each $x_{i,j}$ as a finite linear combination of positive elements $y_{i,j}^{(k)}$, for each of which there are a natural number $m_{i,j}^{(k)}$ and elements $z_{i,j}^{(k,l)}$, $\tilde{z}_{i,j}^{(k,l)}$ such that $y_{i,j}^{(k)} \leq \sum_{l=1}^{m_{i,j}^{(k)}} z_{i,j}^{(k,l)}$ and $z_{i,j}^{(k,l)} \tilde{z}_{i,j}^{(k,l)} = z_{i,j}^{(k,l)}$. Denoting by $e_{i,j}$ the matrix units for $M_n$, it follows that the same is true if we replace the elements $x_{i,j}$,  $y_{i,j}^{(k)}$, $z_{i,j}^{(k,l)}$ and $\tilde{z}_{i,j}^{(k,l)}$ with $x_{i,j} \otimes e_{i,j}$,  $y_{i,j}^{(k)}\otimes e_{i,j}$, $z_{i,j}^{(k,l)} \otimes e_{i,j}$ and $\tilde{z}_{i,j}^{(k,l)} \otimes e_{j,j}$ respectively. It follows that for every $1 \leq i,j \leq n$ $x_{i,j} \otimes e_{i,j}$ belongs to $\Ped (M_n (A))$ and so $x \in \Ped (M_n (A))$. Thus $\Ped (M_n (A)) = M_n (\Ped (A))$. In particular, considering $A$ as a $C^*$-subalgebra of $M_n (A)$ via the embedding $a \mapsto a \otimes e_{1,1}$, it follows that $A \cap \Ped (M_n (A)) = A \cap (M_n (\Ped (A)) = \Ped (A)$.\\
Thus $ \her (a) \subset A \cap M_n (\Ped (A))= \Ped (A)$ and so $a$ belongs to $\Ped (A)$.\\
The proof that $J$ is an ideal goes as in Lemma 3.12 of \cite{purinf} and we omit it.\\
 Since every element in $A$ is a linear combination of positive elements and $J$ is a linear space, we just have to prove that $J_+$ is dense in $A_+$. Let $a \in A_+$, then $\| a - (a-\epsilon)_+\| < \epsilon$ and $[(a-\epsilon)_+ ] \ll [a]$.\\
Hence $J$ is a dense two-sided ideal and has to be equal to $\Ped (A)$. $\Box$\\

\begin{corollary}
\label{cor_ped}
Let $A$ be a $\sigma$-unital stable $C^*$-algebra with almost stable rank $1$ and let $h$ be a strictly positive element for $A$. Then
\[
\begin{split}
\Ped(A) &= \{ a \in A \; : \; \exists b \in A_+ \quad \mbox{ such that } \; [|a|] \ll [b]\}\\
		&=\{ a \in A \; : \; [|a|] \ll [h]\}.
\end{split}
\]
\end{corollary}
\proof 
Since $A$ is stable, by \cite{large} Lemma 2.6, it is large in $A\otimes \mathbb{K}$ in the sense of \cite{large} Definition 2.4. Suppose there are $a \in A$, $b \in M_n (A)_+$ such that $[|a|] \ll [b]$, then there is $\epsilon >0$ such that $[|a|] \leq [(b-\epsilon)_+] \ll [(b-\epsilon /2)_+]$; by largeness, there is $a_{\epsilon /2} \in A_+$ such that $[(b-\epsilon/2)_+] =[a_{\epsilon/2}]$ and thus $[|a|] \ll [a_{\epsilon /2}]$. The result follows from Proposition \ref{ped}. $\Box$\\

\begin{proposition}
\label{ped_f}
Let $A$ be a simple $C^*$-algebra with almost unperforated Cuntz semigroup such that $A\otimes \mathbb{K}$ has almost stable rank $1$ and assume that $A$ admits a unique lower semicontinuous nontrivial $2$-quasitrace $\tau$ (up to scalar multiples).  Let $d_\tau$ be the corresponding functional on $Cu(A)$. Then
\[
\quad \Ped (A) = \{ a \in A \; : \; d_\tau ([|a|]) < \infty \}.
\]
\end{proposition}
\proof 
By Lemma \ref{nontrivial} $d_\tau$ is a nontrivial functional on $Cu (A)$ and is the unique such functional, up to scalars. Let $J' := \{ a \in A \; : \; d_\tau ([|a|]) < \infty \} $.\\
First of all notice that $J'$ is an ideal (again, the proof goes as in Lemma 3.12 of \cite{purinf}, using the fact that $d_\tau$ is a functional on $Cu (A)$) and it is obviously dense since $A$ is simple.\\
By simplicity, every element in $Cu (A)$ is either compact or soft.\\
Let $a \in A$ be such that $0 \neq d_\tau ([| a|])< \infty$; if $[|a|]$ is compact, then $a$ belongs to $\Ped (A)$ by Proposition \ref{ped}, thus we can suppose that $[|a|]$ is soft.\\
Since $0 <d_\tau ([|a|]) <\infty$ and $d_\tau ([|a|]) = \sup_n d_\tau ([ (|a|-1/n)_+])$, there exists $m \in \mathbb{N}$ such that $0 < d_\tau ([ (|a|-1/m)_+]) < \infty$ and so there is a $k \in \mathbb{N}$ such that $k d_\tau ([(|a|-1/m)_+]) \geq d_\tau ([|a|])$. Using comparison and the fact that the only possible functionals on $Cu (A)$ are $\lambda_0$, $\lambda_\infty$ and $d_\tau$ by Lemma \ref{nontrivial}, we see that $[|a|] \leq [(|a| - 1/m)_+ \otimes 1_k] \ll [(|a| - 1/2m)_+ \otimes 1_k]$ and so $a$ belongs to $\Ped (A)$ by Proposition \ref{ped}. The proof is complete. $\Box$\\

Let $A$ be a $C^*$-algebra and $[a] \in Cu(A)$. Then $[a]$ is finite if $[a]+[b]=[a]$ implies $[b]=0$, and $[a]$ is infinite if there is a non-zero element $[b] \in Cu(A)$ such that $[a]+[b]=[a]$ (cfr. \cite{cuntz_t} Section 5.2.2). $[a]$ is properly infinite if $2[a]=[a]$. If $A$ is simple, then $[a] \in Cu(A)$ is infinite if and only if it is properly infinite.

\begin{proposition}
\label{prop2.10}
Let $A$ be a simple $C^*$-algebra with almost unperforated Cuntz semigroup such that $A\otimes \mathbb{K}$ has almost stable rank $1$ and suppose that $A$ admits a unique lower semicontinuous $2$-quasitrace $\tau$ (up to scalar multiples). Then
\[
\Ped (A) = \{ a \in A \;| \; [|a|] \mbox{ is finite }\}.
\]
\end{proposition}
\proof 
Let $a$ be a positive element in $A$. We distinguish different cases.\\
If $[|a|]$ is compact, then $a$ belongs to $\Ped (A)$ by Proposition \ref{ped} and, since $Cu(A)$ admits a non-trivial functional, it cannot contain infinite compact elements by simplicity.\\
Suppose that $[|a|]$ is soft. By simplicity, any infinite element is properly infinite. Thus, if $d_\tau ([|a|]) < \infty$, then $[|a|]$ is finite.\\
If $d_\tau ([|a|])=\infty$, it follows from Proposition \ref{ped_f} that $a$ does not belong to $\Ped (A)$ and by comparison that $|a| \sim |a| \oplus |a|$; hence $[|a|]$ is properly infinite. $\Box$\\

In \cite{cuntz_hm} it is defined a notion of Cuntz subequivalence and compact containment for countably generated Hilbert modules of a $C^*$-algebra and is proved that the induced partially ordered semigroup coincides with the Cuntz semigroup, in the sense that the two constructions are naturally equivalent (\cite{cuntz_hm} Section 6). Given two countably generated Hilbert modules $X$ and $Y$ over a $C^*$-algebra, we write $X \precsim Y$ if $X$ is Cuntz subequivalent to $Y$, and $X \Subset Y$ if $X$ is compactly contained in $Y$; we refer to \cite{cuntz_hm} Section 1 for the definitions.

\begin{theorem}
\label{dyc2}
Let $A$ be a $\sigma$-unital stable simple $C^*$-algebra with almost unperforated Cuntz semigroup, with almost stable rank $1$ and admitting a unique lower semicontinuous nontrivial $2$-quasitrace $\tau$. Let $x$ be a positive element in $A$, then we have
\[
[x] \; \mbox{ is infinite } \Leftrightarrow \her(x) \; \mbox{ is stable }  \Leftrightarrow E_x  \simeq A
				\]
and
\[
[x] \; \mbox{ is finite } \Leftrightarrow \her(x) \; \mbox{ is algebraically simple } \Leftrightarrow E_x \Subset A.
\]
Furthermore, if $y$ is a positive element in $A$ such that $[x]<[y]$, then $[x] \ll [y]$ and this happens if and only if there exist elements $z \in A$ and $e$ self adjoint in $\her(y)$ such that $\her(x) = \her (z^* z)$, $\her(zz^*) \subset \her(y)$, $ea=a$ for every $a \in \her(zz^*)$. This condition is also equivalent to the existence of a Hilbert module $E$ over $A$ such that $E_x \simeq E \Subset E_y$.
\end{theorem}
\proof
If $A$ is simple, the Cuntz class of a positive element $x \in A$ is infinite if and only if it is properly infinite.
Now, if $[x]$ is properly infinite, then it defines the same Cuntz class of any strictly positive element for $A$; by \cite{cuntz_hm} Theorem 3 and \cite{open_proj} Proposition 4.3 it follows that, since $A$ has almost stable rank $1$, $\her (x) \simeq A$ and $E_x \simeq A$ as countably generated Hilbert modules.\\
Conversely, the Cuntz class of a strictly positive element in a stable $C^*$-algebra is properly infinite, as well as the generator of a Hilbert module that is isomorphic to $A$ and so the characterization of infinite elements in the Cuntz semigroup in terms of hereditary subalgebras and Hilbert modules is complete.\\

By Proposition \ref{prop2.10} $[x]$ is finite if and only if $x$ belongs to $\Ped (A)$ and in this case $\her (x)$ is algebraically simple by \cite{pedersen} Proposition 5.6.2. Note now that if $\her (x)$ is algebraically simple and $x$ does not belong to $\Ped (A)$, then by Proposition \ref{ped_f} $d_\tau ([x])=\infty$, then it follows that $\her(x)$ only admits the trivial lower semicontinuous $2$-quasitrace $\tau_\infty$, contradicting the hypothesis. Hence, under our assumptions, $\her (x)$ is algebraically simple if and only if $[x]$ is finite.\\
By Corollary \ref{cor_ped}, $x$ belongs to $\Ped(A)$ if and only if $[x] \ll [h]$ for any strictly positive element $h$ for $A$. From the equivalence between compact containment for classes of Hilbert modules and the way below relation (see \cite{cuntz_hm} Theorem 1) it follows that $[x] \ll [h]$ if and only if $E_x \Subset A$ up to isomorphism of Hilbert modules.\\

Suppose now that there is a positive element $y \in A$ such that $[x]<[y]$. If either $[x]$ or $[y]$ is compact, then $[x] \ll [y]$, hence we can suppose $[x]$ and $[y]$ are soft. Note that $d_\tau ([(y-\epsilon)_+]) < d_\tau ([y])$ for every $\epsilon >0$, otherwise $[y]$ would be compact by comparison and so, since $d_\tau ([y])= \sup_{\epsilon >0} d_\tau ([(y-\epsilon)_+])$, there exists $\eta>0$ such that $d_\tau ([x]) \leq d_\tau ([(y-\eta)_+])$, entailing $[x] \ll [y]$. From the equivalence between compact containment in the concrete and in the abstract sense proved in \cite{cuntz_hm} Theorem 1, this is the same as the existence of a Hilbert module $E$ such that $E_x \precsim E \Subset E_y$. Thus, in the case $A$ has almost stable rank $1$, it follows from Theorem 3 of \cite{cuntz_hm} that there exists a Hilbert module $E'$ such that $E_x \simeq E' \subset E \Subset E_y$, hence in particular $E_x \simeq E' \Subset E_y$. Denote by $\Phi : E_x \rightarrow E'$ an isomorphism of Hilbert modules; since $A$ is $\sigma$-unital, $x$ belongs to $E_x$ and $E' = \overline{\Phi (x) A}$, with $\Phi (x) \in A$ since $E' \subset E_y$; it follows then from Proposition 4.3 of \cite{open_proj} that there is an element $z \in A$ such that $\her (x)=\her (z^* z)$ and $\her (\Phi (x)) = \her (zz^*)$. If now $e \in \mathbb{K}(E_y)= \her (y)$ is a self adjoint element such that $e|_{E_{\Phi(x)}} = \id|_{E_{\Phi(x)}}$, it follows that $ea=a$ for every $a \in \her (\Phi (x))$. Suppose now that there are an element $z \in A$ and $e$ selfadjoint in $\her (y)$ such that $\her (x) = \her (z^* z)$, $\her (zz^*) \subset \her (y)$ and $ea=a$ for every $a \in \her (zz^*)$. Then again by Proposition 4.3 of \cite{open_proj} it follows that $E_x \simeq E_{zz^*}$ and since $zz^*$ belongs to $\her (zz^*)$, it follows that $ezz^* = zz^*$, which implies $e|_{E_{zz^*}} =  \id|_{E_{zz^*}}$ and hence $E_x \simeq E_{zz^*} \Subset E_y$. In this case, it follows again by \cite{cuntz_hm} Theorem 1 that $[x]\ll [y]$. $\Box$\\

\begin{remark}
\label{exaf}
There are examples of non-algebraically simple, non-stable simple $\sigma$-unital $C^*$-algebras with stable rank $1$, admitting non-trivial lower semicontinuous $2$-quasitraces and whose Cuntz semigroup is almost unperforated. They were constructed by Blackadar in \cite{af_dyc}; namely, in Corollary 4.5 it is proved that for any simple unital $AF$-algebra $A$ and any $F_\sigma$ subset of its space of bounded traces, it is possible to find an hereditary $AF$-subalgebra $B$ of $A\otimes \mathbb{K}$ containing $A$ such that the bounded traces on $B$ are exactly the extensions coming from the $F_\sigma$ set chosen. It is clear that this construction does not work if $A$ has just one trace.
\end{remark}

\begin{example}
The following example is well known. The $C^*$-algebra $\mathbb{K}$ of compact operators on a separable Hilbert space has stable rank $1$ and a unique nontrivial trace, up to scalars. Its Cuntz semigroup is $Cu (\mathbb{K})= \bar{\mathbb{N}}$, which is almost unperforated. By Theorem \ref{dyc2} any hereditary $C^*$-subalgebra of $\mathbb{K}$ is either algebraically simple (actually unital) or isomorphic to $\mathbb{K}$.
\end{example}

\begin{example}
The Jiang-Su algebra $\mathcal{Z}$ (see \cite{js}) has stable rank $1$, a unique trace and almost unperforated Cuntz semigroup equal to $Cu (\mathcal{Z})= \mathbb{N} \sqcup (0,\infty]$. The hereditary $C^*$-subalgebras of $\mathcal{Z} \otimes \mathbb{K}$ generated by compact elements (corresponding to $\mathbb{N} \subset Cu (A)$) are unital; the ones corresponding to elements in $(0,\infty)$ are nonunital and algebraically simple.
\end{example}

\begin{example}
The $C^*$-algebra $\mathcal{W}$ introduced by Jacelon in \cite{razak} has stable rank $1$, a unique trace and its Cuntz semigroup is $Cu (\mathcal{W})= [0,\infty]$, which is almost unperforated (one could also just notice that $\mathcal{W}$ is $\mathcal{Z}$-stable). Thus any $\sigma$-unital $C^*$-algebra that is Morita-equivalent to $\mathcal{W}$ is either algebraically simple (in correspondence of elements in $(0,\infty) \subset Cu (\mathcal{W})$) or stable.
\end{example}

\begin{example}
If $X$ is a compact manifold with finite covering dimension foliated by a free minimal action of $\mathbb{R}$ admitting a unique invariant transverse measure, then $C(X) \rtimes \mathbb{R}$ satisfies the hypothesis of Theorem \ref{dyc2} in the case it is projectionless. This is a consequence of \cite{ms} Theorem 6.30, \cite{rokflow} Corollary 9.2 and \cite{zstable_projless} Corollary 3.2.
\end{example}

\subsection*{Acknowledgements}
This research was partially funded by I.N.d.A.M. (grant n. U-UFMBAZ-2018-000926) and the University of Tor vergata (grant n. 669240 CUP: E52I15000700002). The author thanks Prof. Longo for the hospitality at the Universit\`a degli Studi di Roma Tor Vergata for the period of this research. Special thanks go to Andr\' e Schemaitat and Dominic Enders for many interesting discussions. In particular, Proposition \ref{ped} is the result of fruitful discussions with the colleague and friend  Andr\' e Schemaitat.

\end{document}